\documentclass[leqno,12pt]{article}
\usepackage{epic}
\usepackage{ amsmath,amstext,amsthm,amssymb}
\numberwithin{equation}{section}

\begingroup   

\allowdisplaybreaks
\swapnumbers
\theoremstyle{plain}

\def\bthm#1.#2 #3\ethm{ 
\begin{\ifx#1ttheorem\fi%
\ifx#1llemma\fi%
\ifx#1ccorollary\fi%
\ifx#1pproposition\fi%
\ifx#1ddefinition\fi} \label{#1.#2} 
#3 \end{\ifx#1ttheorem\fi%
\ifx#1llemma\fi%
\ifx#1ccorollary\fi%
\ifx#1pproposition\fi%
\ifx#1ddefinition\fi}} 
 
\def\t#1/{theorem~\ref{t#1}}   \def\T#1/{Theorem~\ref{t#1}} 
\def\c#1/{corollary~\ref{c#1}}   \def\C#1/{Corollary~\ref{c#1}} 
\def\l#1/{lemma~\ref{l#1}}        \def\L#1/{Lemma~\ref{l#1}} 
\def\s#1/{section~\ref{s#1}} 
\def\e#1/{(\ref{e#1})} 
\def\d#1/{definition~\ref{d#1}} 
\def\f#1/{figure~\ref{f#1}} 
 
\def\Label #1 {\label{#1}}

\def\norm#1.#2.{\lVert#1\rVert_{#2}} 
\def\Norm#1.#2.{\bigl\lVert#1\bigr\rVert_{#2}} 
\def\NOrm#1.#2.{\Bigl\lVert#1\Bigr\rVert_{#2}} 
\def\NORm#1.#2.{\biggl\lVert#1\biggr\rVert_{#2}} 
\def\NORM#1.#2.{\Biggl\lVert#1\Biggr\rVert_{#2}} 
 
\def\ip#1,#2.{\langle #1,#2\rangle} 
\def\Ip#1,#2.{\bigl\langle#1,#2\bigr\rangle} 
\def\IP#1,#2.{\Bigl\langle#1,#2\Bigr\rangle}

\def\abs#1{\lvert#1\rvert} 
 
\def\ABs#1{\Bigl\lvert#1\Bigr\rvert} 
\def\ABS#1{\biggl\lvert#1\biggr\rvert}

\newcommand{\zc}{\ensuremath{\psi}} 
 
\newcommand{\zd}{\ensuremath{\delta}}

\newcommand{\zve}{\ensuremath{\varepsilon}} 
\newcommand{\zvf}{\ensuremath{\varphi}} 
\newcommand{\zf}{\ensuremath{\phi}} 
\newcommand{\zF}{\ensuremath{\Phi}} 
\newcommand{\zg}{\ensuremath{\gamma}} 
\newcommand{\zG}{\ensuremath{\Gamma}} 
 
\newcommand{\zI}{\ensuremath{\infty}}

\newcommand{\zm}{\ensuremath{\mu}}

\newcommand{\zt}{\ensuremath{\tau}} 
\newcommand{\zw}{\ensuremath{\omega}} 
 
\newcommand{\zx}{\ensuremath{\xi}} 
 
\newcommand{\zs}{\ensuremath{\sigma}}

\newcommand{\zp}{\ensuremath{\pi}} 
\newcommand{\zq}{\ensuremath{\chi}}

\def\z#1#2{\ifcase#1 {\mathcal {#2}}  
\or {\mathbf{#2}}                    
\or  {\boldsymbol{#2}}                  
\or{\widetilde{#2}}                   
\or {\acute{#2}}
\or {#2 }		
\or  {\mathbf{#2}_{\text{\rm fat}}}	
\or   {#2}_{\text{\rm fat}} 		
\or  {#2}_{\text{\rm thin}}            
\or {\mathbf{#2}}_{\text{\rm thin}}\fi}      
 
\def\ZR{\ensuremath{\mathbb R}} 
 
\def\ZZ{\ensuremath{\mathbb Z}}

\def\seq{\preceq }   
\def\mid{\,:\,}

\def\md#1#2\emd{\ifx0#1 
\begin{equation*} #2 \end{equation*}\fi  
\ifx1#1\begin{equation}#2\end{equation}\fi   
\ifx2#1\begin{align*}#2\end{align*}\fi   
\ifx3#1\begin{align}#2\end{align}\fi    
\ifx4#1\begin{gather*}#2\end{gather*}\fi  
\ifx5#1\begin{gather}#2\end{gather}\fi   
\ifx6#1\begin{multline*}#2\end{multline*}\fi  
\ifx7#1\begin{multline}#2\end{mutline}\fi  
}

\begin{document}

\title{   Carleson's Theorem
\\ with Quadratic Phase Functions}

 \author{Michael T. Lacey\thanks{This work has been 
  supported by an NSF grant, DMS--9706884.}\\
  Georgia Institute of Technology }

\maketitle

\def\ex#1]{e(#1)}


\section{The Main Result}
Consider the maximal operator 
\md0
\z0C_d f(x)=\sup_{\text{deg}(p)=d}\ABS{ \text{p.v.} 
\int f(x-y)\ex p(y)] \frac{dy}y} 
\emd
in which $d$ is an integer, $p$ is a polynomial of degree $d$, $\ex 
u]:=e^{\zp i u}$, $f$ is a
Schwarz function and the integral is understood in the principal value 
sense.  
This definition is motivated principally by the case 
$d=1$.
$\z0C_1f$ controls the maximal partial Fourier 
integrals of $f$ and it  
extends to a bounded map from $L^p$ into itself for $1<p<\zI$.  The 
critical contribution here is 
L.\thinspace Carleson's proof \cite{c} of the boundedness of $\z0C_1$ 
from $L^2$ into weak--$L^2$.
The $L^p$ version was established by R.~Hunt \cite{h}.  Also see 
\cite{f,ltcar}.  

It is natural to ask if the same results hold for larger values of $d$.  
Indeed, it does for 
the case of $d=2$ and this is the main result of our paper.

\bthm  t.main   $\z0C_2$ extends to a bounded map from $L^p$
 into   itself for all $1<p<\zI$.
\ethm

To prove the theorem, it suffices to show that $\z0C_2$  maps $L^2$ into weak
$L^2$ as 
  our proof can be modified to treat all $1<p<\zI$, and we briefly indicate how 
to do this in the next section.  

E.\thinspace M.\thinspace Stein \cite{st} has observed that the Fourier transform
of the distribution
$\ex y^2]/y$ has an easily calculable form, one that permits analysis
of the maximal
operator formed from dilations of this distribution.  We relie very much
on this observation.
As well, we now have   a much richer understanding of Carleson's theorem as presented in
papers of Fefferman, Lacey and Thiele  and Thiele
\cite{f,ltcar,t} and some related issues \cite{lt1,lt2,lmax}.   We invoke some of these
elements to provide a proof of our main theorem.    Stein's argument
and the overview of our proof
are laid out in the next section.  The main inequality described in
that section requires a careful
analysis in time and frequency variables, that being carried out in the
remaining sections of the paper.

\smallskip

The elegant results of K.~Oskolkov \cite{o} are  of the same genre as ours.

\smallskip

T.~Tao and J.~Wright informed me of this problem.  G.~Mockenhoupt brought
Stein's article \cite{st} to my attention.   Part of this work
was completed at the Centre for Mathematics and its Applications
at the Australian National University.  I am endebited to all.

\medskip


Notations:  The Fourier transform is taken to be $\widehat f(\zx):=\int
f(x)\ex -2x\zx]\; dx$.
The norm of an operator $T$ from $L^p$ into $L^p$ is written as $\norm
T.p\to p.$ with a
corresponding notation for the weak type norm.   By $A\seq{}B$ we mean
that there is an absolute
constant $K$ so that $A\le{}KB$.  By $A \simeq{}B$ we mean $A\seq{}B$
and $B\seq{}A$.  $c(J)$ is the center of the interval of $J$.
$\zq_J(x)=(1+\text{dist}(x,J)\abs{J}^{-1})^{-1}$.


\section{The Overview of the Proof}
The supremum we 
 wish to bound admits a description as a supremum over dilations, for which 
there are a wealth of techniques to use, and a supremum over modulations in 
frequency, which is the domain of Carleson's theorem.  
It is useful to formalize these aspects with a couple of definitions.

A distribution $K(y)$  determines two maximal functions of interest to us.
They are
\md2
\z0D[K]f(x):={}&\sup_{a>0}\ABS{\int f(x-y)a^{-1}K(a^{-1}y)\;dy} ,
\\
\z0C[K]f(x):={}&\sup_{b\in\ZR}\abs{\z0D[K](\ex b\cdot]f)(x)}
\\{}=& \sup_{b\in\ZR}\sup_{a>0}\ABS{\int f(x-y)\ex
by]a^{-1}K(a^{-1}y)\;dy} .
\emd
Thus if $K(y)=y^{-1}1_{\{0<\abs y\le1\}}$, $\z0D[K]$ is the maximal
truncations of the
Hilbert transform and $\z0C[K]$ is Carleson's maximal operator.   Set
$K(y):=\ex y^2/4]/y$.  To prove
our theorem we show that $\z0C[K]$ maps $L^2$ into $L^{2,\zI}$.

We recall Stein's argument \cite{st} that $\z0D[K]$ maps $L^2$ into $L^2$.  The
Fourier transform of $K$ is a
smooth odd function  satisfying
\md1\Label e.k1
\widehat K(\zx)=c_0+\ex \zx^2]\{ c_1/\zx+c_2/\zx^2+\cdots\}\qquad \text{as
$\zx\to\zI$}
\emd
 for some choice of constants $c_j$, $j\ge1$.  Indeed $\partial \widehat
K(\zx)={} \int \ex x^2/4-x\zx]\; dx={} c\ex \zx^2]$.
Moreover $\widehat K$ is odd as $K$ is odd hence
\md0
c^{-1}\widehat K(\zx)=\int_0^\zx\ex y^2]\; dy={}
\int_0^\zI\ex y^2]\; dy-\int_\zx^\zI\ex y^2]\; dy.
\emd
And the assertation follows as
\md0
\int_\zx^\zI\ex y^2]\; dy=\ex \zx^2]\{ c_1'/\zx+c_2'/\zx+\cdots\}\qquad
\text{as $\zx\to\zI$}
\emd

\smallskip

With \e.k1/ established we can write $\widehat K=\widehat H+\sum_{j=1}^\zI
2^{-j}\widehat m_j$ where $\widehat H$ is
 smooth odd and equals $c_0\text{sign}(\zx)$ for $\abs\zx>2$.  The
multipliers $\widehat m_j$ are of the form
\md1 \Label e.mjdef   \widehat m_j(\zx)=\zc(2^{-j}\zx)\ex \zx^2]
\emd
where $\zc$ is a $C^\zI$ function with support in
$\frac12\le\abs\zx\le2$.   These multipliers are our main concern.  

Now, the term arising from $H$ is governed by the Hilbert transform.
In particular $\z0D[H]$ is bounded from $L^2$ into $L^2$ as this is the
maximal truncations of the Hilbert transform.  In a like manner we have
the estimate $\norm \z0C[H].2\to2,\zI.\seq1$ by Carleson's Theorem \cite{c}.

Stein has shown that $\norm \z0D[m_j].2\to2.\seq2^{j/2}$, which then
completes the proof of the
bound on $\z0D[K]$.
To prove our main theorem, we demonstrate that
\md1\Label e.C
\norm \z0C[m_j].2\to2,\zI.\seq{}2^{\zg j} \qquad\text{for some
$0<\zg<1$.}
\emd
In fact $\zg=8/9$ will work.  [By optimizing our argument we could
establish this estimate for any $\zg>1/2$.]

\smallskip
Stein's argument is crucial to our own and so we recall it here.  To
bound $\norm \z0D[m_j].2\to2.$ it suffices to prove that
\md1\Label e.mj
\NOrm \sup_{1\le{}a\le2}\ABs{\int m_{j,a}(y)f(x-y)\;
dy}.2.\seq{}2^{j/2}\norm f.2.
\emd
where $m_{j,a}(y):=a^{-1}m_j(a^{-1}y)$.  And to this end the method of
$TT^*$ is invoked.  Observe that
\md5\Label e.mja  \sup_{1\le a,b\le 2}
\abs{m_{j,a}*\overline{m_{j,b}(y)}}\seq{}\zF_j(y),\quad y\in\ZR,\
\\
 \Label e.fj
\text{where}\quad \zF_j(y)= \begin{cases}
2^{j/2}\abs{y}^{-1/2} & \abs{y}\le c2^j  \\
2^j(1+2^j\abs{y})^{-2} &  \abs y \ge c2^j    \end{cases}  .
\emd
We can take $c=16$.  Note that $\zF_j$ is non decreasing and $\int \zF_j\;
dy\seq2^j$ which proves \e.mj/.

After taking dilation into account, \e.mja/ amounts to the estimate
\md1\Label e.Fj
\sup_{0\le{}b\le{}4}\ABs{  \int \tilde\zc(2^{-j}\zx)\ex b\zx^2+2\zx y]\;
d\zx}\seq{}
\zF_j(y).
\emd
Here $\tilde \zc$ is another Schwarz function with support in
$\frac12\le{}\abs\zx\le2$.  Set $p(\zx)=b\zx^2+2\zx y$.  If $\abs y\ge{}cb2^j$
observe that the derivative of $p$ with respect to $\zx$ exceeds $\overline c
\abs y$ on the support of $\tilde\zc(2^{-j}\zx)$.   Thus repeated integration by
parts will prove the estimate.  If $\abs y\le{}cb2^j$
we can use   the van der
Corput second derivative test.
It provides the estimate of the integral as $\seq{}b^{-1/2}\seq
2^{j/2}\abs{y}^{-1/2}$.  Thus the inequality holds.

\medskip

 The reminder of the paper is devoted to a
proof of \e.C/.

To do so we use the time frequency analysis of
Lacey--Thiele \cite{ltcar} with some further ideas drawn from
Fefferman and Thiele \cite{f,t}. A central conceptual problem
arises from the fact that $\widehat m_j$ is supported in an
interval of length $2^j$ but $m_j$ has (approximate) spatial
support in an interval of the same length. That is classical
Fourier uncertainty is not observed. Treating this issue is
probably the main novelty of this paper.

For our subsequent use observe these points.  First in the definition
of \e.mjdef/ we can assume that $\zc$ is supported in
$1-\frac1{40}\le\zx\le1+\frac 1{40}$ (as opposed to $\frac12\le\abs \zx\le2$).

Second \e.Fj/ implies that $\int\abs{m_j}\; dy\seq2^j$ hence $\norm
\z0D[m_j].\zI\to\zI.\seq2^j$ so that by interpolation
\md1  \Label   e.mjp
\norm \z0D[m_j].p\to p.\seq2^{j(1-1/p)},\qquad 2<p<\zI.
\emd
Thus $\z0D[K]$ maps $L^p$ into itself for $2<p<\zI $. 

\smallskip 
Indeed, this estimate holds for $1<p<\zI$.  In fact, we have the estimate 
$\norm \z0D[m_j].\zI\to\zI.\seq2^j$, with the same estimate holding at $p=1$ as 
well. These estimates require no cancellation, and so hold for $\z0C[m_j]$ as 
well.  Thus to prove our main theorem and 
 in light of the extension of Hunt of Carleson's theorem, it 
suffices to provide the bound we have claimed for $\z0C[m_j]$ on $L^2$.  

\smallskip

Third there is a sharper form of Stein's observation.  Namely the
operator
\md1\Label  e.dj
D_jf(x):=\sup_{1\le a\le2}\sup_{\abs N\le2^j}
\ABs{\int a^{-1}m_j(a^{-1}(x-y))\ex Ny]f(y)\;dy}
\emd
maps $L^2$ into $L^2$ with norm bounded by ${}\seq2^{j/2}$.  Employing
the same arguments as above, this amounts to the estimate
\md4
\sup_{0\le b\le4}\sup_{\abs N\le 2^{j+2}}
\ABs{\int \tilde\zc(2^{-j}\zx)\ex b\zx^2+(N+y)\zx]\;dy}\seq\zG_j(y)
\\
\text{where}\quad\zG_j(y)=  \begin{cases}
2^{j/2}\abs{y}^{-1/2} & \abs y\le{}c2^j
\\
 2^j(1+2^j\abs y)^{-n} &  \abs y\ge c2^j  \end{cases}
\emd
In this definition, $n$ is an arbitrary positive integer and $c=32$.
Details are a modification
of the earlier argument.  In fact we have $\norm
D_j.p.\seq{}2^{j(1-1/p)}$ for $2\le{}p<\zI$.  We
shall have recourse to this below.

Fourth in proving the estimate \e.mj/ we follow the approach of
Kolomogorov and  Silverstoff, as Fefferman \cite{f} has demonstrated that this is a
powerful technique in issues related to Carleson's theorem.   We show that there is a $0<\zg<1$ so
that for all $j$, measurable functions $N\mid \ZR\to\ZR$, $\ell\mid
\ZR\to\ZZ$ and $a\mid \ZR\to [1,1+\frac1{40}]$,
\md1\Label e.2do
\NOrm  \int m_{j,a(x)2^{\ell(x)}}(x-y)f(y)\ex N(x)y]\; dy
.2,\zI.\seq{}2^{\zg j}\norm f.2. .
\emd
We will do this with $\zg=8/9$.   This inequality is sufficient for our
purposes.


\section{The Discrete Operator}




\def\thin #1]{ \text{\rm thin}(#1)}

Let ${\mathbf D}$ be a collection of dyadic intervals in the real line.   Let
$\z6P$ be the
set of rectangles $\z5s=I_{\z5s}\times \zw_{\z5s}\in{\mathbf D}\times {\mathbf D}$
   which have area $\abs{I_{\z5s}} \abs {\zw_{\z5s}}=2^{2j}$. We call these ``fat
 tiles" and we generically write $s,s',s''$ for fat tiles. Let
$\zw_{\z5s1}$
($\zw_{\z5s2}$) be the left (right) half of $\zw_{\z5s}$.  This
definition is chosen in accordance with the frequency and spatial localizations
of the kernel $m_j$, its dilates and modulations.

Let $\z9P$ be the set of rectangles $s\in{\mathbf D}\times {\mathbf D}$ of area 1.
We call these ``thin tiles" and we generically write $\zs,\zs',\zs''$ for thin
tiles.       Set $\thin s]:=\{\zs\in\z9P\mid
I_\zs=I_{\z5s},\ \zw_\zs\subset\frac34\zw_{\z5s 2}\}$.     For $\zs\in\thin \z5s]$, set $
\zw_{\zs j}:=\zw_{\z5sj}$ for $j=1,2$.  See figure 1.  [Actually,
``thin tiles" obey classical Fourier uncertainty and so are thin
only in contrast to fat tiles.]


\begin{figure}
\begin{center}
\setlength{\unitlength}{0.00023333in}
\begingroup\makeatletter\ifx\SetFigFont\undefined%
\gdef\SetFigFont#1#2#3#4#5{%
  \reset@font\fontsize{#1}{#2pt}%
  \fontfamily{#3}\fontseries{#4}\fontshape{#5}%
  \selectfont}%
\fi\endgroup%
{\renewcommand{\dashlinestretch}{30}
\begin{picture}(10158,7473)(0,-10)
\thicklines
\drawline(750,3975)(10125,3975)
\drawline(825,6600)(10050,6600)(10050,6225)
	(825,6225)(825,6600)
\drawline(825,6225)(10050,6225)(10050,5850)
	(825,5850)(825,6225)
\drawline(825,5850)(10050,5850)(10050,5475)
	(825,5475)(825,5850)
\drawline(825,5475)(10050,5475)(10050,5100)
	(825,5100)(825,5475)
\drawline(825,5100)(10050,5100)(10050,4725)
	(825,4725)(825,5100)
\drawline(825,7425)(10125,7425)(10125,525)
	(825,525)(825,7425)
\put(-370,2475){\makebox(0,0)[lb]{\smash{{{\SetFigFont{10}{12}{\rmdefault}{\mddefault}{\updefault}$\zw_{s1}$}}}}}
\put(-305,5775){\makebox(0,0)[lb]{\smash{{{\SetFigFont{10}{12}{\rmdefault}{\mddefault}{\updefault}$\zw_{s2}$}}}}}
\put(5100,-200){\makebox(0,0)[lb]{\smash{{{\SetFigFont{10}{12}{\rmdefault}{\mddefault}{\updefault}$I_s$ }}}}}
\end{picture}
}
\end{center}
\Label   f.fatthin
\caption{\small A fat tile $s$ and thin tiles $\zs\in\thin s]$.}
\end{figure}
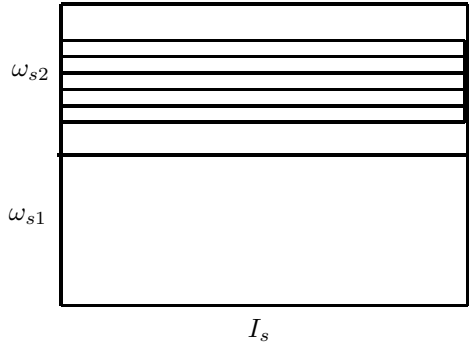

Fix a Schwarz function $\zf$ with
$1_{[-\frac1{80},\frac1{80}]}\le\widehat\zf\le1_{[
-\frac1{40},\frac1{40}]}$.  For a rectangle $\zs=I_\zs\times \zw_\zs$ of area 1
(not necessarily a thin tile) define
\md0
\zf_\zs(x)=\frac{\ex
c(\zw_\zs)x]}{\sqrt{\abs{I_\zs}}}\zf\Bigl(\frac{x-c(I_\zs)}{\abs{I_\zs}}
\Bigr).
\emd
In this display and throughout, $c(J)$ is the center of the interval
$J$.

Fix the data $j\ge1$, $f\in L^2$ of norm one, functions $N$, $\ell$ and
$a$ as in \e.2do/.
For $\z5s\in\z6P$, $\zs\in\thin \z5s]$ and integer $l$ with
$2^{-l}=\abs{\zw_{\z5s2}}$.  Define
\md4
E(\z5s)=E(\zs):=\{x\in\ZR\mid N(x)\in\zw_{\z5s1},\  \ell(x)=l\},
\\
\zvf_\zs(x)=1_{E(\zs)}(x)\int\zf_\zs(y)\ex
-2N(x)y]a(x)^{-1}2^{-l-j}m_j(a(x)^{-1}2^{-l-j}(x-y))\; dy ,
\\
M_jf(x)=\sum_{\z5s\in\z6P}\sum_{\zs\in\thin \z5s]}\ip f,\zf_\zs.\zvf_\zs(x).
\emd

A principal motivation for these definitions is the proof of \l.tree/
below.  At this point we simply observe that the support of the integral in the
definition of  $\zvf_\zs$
is in $E(\zs)$.  $\widehat m_j$ is supported in a small neighborhood of $2^j$ so
that the second function in the last integral has frequency support in a small
interval around $2^{-l}$.  $\widehat \zf_\zs$ is supported in a small interval
around $c(\zw_{\zs})$ with $\zw_\zs\subset\frac34\zw_{\zs2}$.  So $N(x)$ must be in
$\zw_{s1}$ in order for the integral to be non zero.

We claim that the following inequality is sufficient for \e.2do/.
\md1  \Label e.2DO
\norm M_jf.2,\zI.\seq2^{\zg j},\qquad \zg=8/9.
\emd
In the proof of this inequality, we 
 only consider sums
over finite subsets $\z6S\subset\z6P$.
We fix data $f\in L^2$ of norm one and the functions $N$, $\ell$ and
$a$.  Let $M_j$ be the sum restricted to this new smaller class of tiles.
Then, by dilation invariance, \e.2DO/ is implied by  this
inequality. \md1 \Label e.d2o
\abs{\{ M_jf>1\}}\seq2^{2\zg j},\qquad \zg=8/9,
\emd
the inequality holding for all functions $f $ of norm one.

\begin{proof}[Proof of sufficency of \e.2DO/
.] 
A convexity argument can be used to show that \e.2DO/
implies the inequality \e.C/.  
Indeed arguments like this have been used many times in related papers, for 
instance \cite{f,ltcar}.

Let us give the convexity argument in an elemental form.  For our subsequent 
use, let us define translation and modulation operators by $Tr_t f=f(x-t)$ and 
$Mod_tf(x)=e^{-ixt}f(x)$ for $t\in\ZR$.  Observe that the sum
\md0
\sum_{n\in\ZZ}\ip f, Tr_n \zf. Tr_n\zf
\emd
could bewritten as a sum over tiles.  
More importantly, 
\md2
\int_0^1 \sum_{n\in\ZZ}\ip f, Tr_{n+t} \zf. Tr_{n+t}\zf
\;dt={} &
 \int_{-\zI}^\zI \ip f, Tr_t \zf. Tr_t\zf \; dt
 \\
 {}={}&\zc*f
 \emd
 where $\zc(x)=\int\overline{\zf(y)}\zf(x+y)\; dy$.  
 Recall that we specified $\zf$ to be a Schwartz function with $
 1_{[-\frac1{80},\frac1{80}]}\le\widehat\zf\le1_{[
-\frac1{40},\frac1{40}]}$, so that $\zc$ satisfies a similar set of 
inequalities. 

Elaborating on this theme, observe that this sum 
\md1 \Label e.modt \sum_{m,n\in\ZZ} \ip f, Mod_mTr_n\zf.Mod_mTr_n\zf
\emd
could be written as a sum over tiles.  
Define 
\md2
Af={}&
\int_0^1 \int_0^1\sum_{m,n\in\ZZ}
\ip f, Mod_{m+\zt}Tr_{n+t}\zf.Mod_{m+\zt}Tr_{n+t}\zf   \;d\zt dt
\\{}={}& 
  \int_{-\zI}^{-\zI} \int_{-\zI}^{-\zI}
\ip f, Mod_{  \zt}Tr_{  t}\zf.Mod_{  \zt}Tr_{  t}\zf   \;d\zt dt
\emd
This is a multiple of the identity, as is easy to see.  

By periodicity, $Af$ is also equal to 
\md1 \Label e.period  
Af=\lim_{\z2\zt\to\zI}\lim_{\z1t\to\zI} (\z2\zt\z1t)^{-1}
  \int_{0   }^{\z1t} \int_{  0 }^{\z2 \zt
   }\sum_{m,n\in\ZZ}
\ip f, Mod_{m+\zt}Tr_{n+t}\zf.Mod_{m+\zt}Tr_{n+t}\zf   \;d\zt dt
\emd
This concludes our general remarks on the use of convexity.  

\medskip

Let us turn to the operator $M_j$.  Define, for an integer $l$ 
\md0
P_{j,l}f=\sum_{\substack{ s\in \z6P\\ \abs{\zw_{s2}}
=2^{-l} }}\sum_{\zs\in\thin s  ]} \ip f,\zf_\zs.\zf_\zs ,
\emd
and observe that this
sum is similar to \e.modt/.  We may average these operators 
over modulations and translations to obtain a multiple of the identity.  This 
can be done in a way that is independent of $l\in\ZZ$ and essentially independent 
of $j\ge1$. 
We shall return to this point momentarily.

To make the connection with our operator  $M_j$ 
more directly, observe that with the notation used in the definition of $M_j$, 
\md0
M_j f(x)=( Mod_{2N(x)}Dil_{2^{l(x)-j}a(x)}m_j)*P_{j,\ell(x)}f(x)
\emd
where $Dil_\zd {}g(x)=\zd^{-1}g(x\zd^{-1})$.  

Thus the main point is that we can recover the identity operator from $P_{j,l}$ 
in a way that is independent of $l$ and $j$ and does not effect the assumed 
inequality \e.2DO/.  

But certainly translation and modulation do not efect the distributional 
inequality.  And, we can obtain the identity operator from the $P_{j,l}$ in this 
way.  Recall that the tiles depend upon choices of dyadic grids $\z1D$ and 
$\z1D'$.  A translation of $\z1D$ ($\z1D'$) corresponds to an application of 
$Tr_t$ ($Mod_\zt$) to the functions $\zf_\zs$.  Thus the assumed inequality 
applies   to any $M_j$ obtained from translations of either grid.  
Finally, the periodicity 
property \e.period/ shows that the identity operator can be obtained in a 
way that is independent of $l$.  This completes the proof. 

\end{proof}


\section{Trees and size}
The principle definitions and lemmas are stated in this section.  We
show how they prove \e.d2o/ and prove the Lemmas in the following
section.    We begin with requisite definitions.

For $\z5s,\z5s'\in\z6P$ say that $\z5s<\z5s'$ iff $I_{\z5s}\subset
I_{\z5s'}$ and
$\zw_{\z5s}\supset \zw_{\z5s'}$.  Say that $\z7T\subset \z6P$ is a {\em
tree} if there is a $I_{\z7T}\times
\zw_{\z7T}\in\z6P$ with $\z5s<I_{\z7T}\times \zw_{\z7T}$ for all
$\z5s\in\z7T$.

A subset $\z8T\subset\z9P$ is a {\em tree} if it is a subset of
$\thin \z7T]$ for some tree $\z7T\subset\z6P$.   We   denote
 the top of
the tree by $I_{\z8T}\times \zw_{\z8T}$.
A tree $\z8T$ is a $1$--tree ($2$--tree) iff  for all $\zs,\zs'\in\z8T$
either $\zw_{\zs}=\zw_{\zs'}$ or $\zw_{\zs}\cap\zw_{\zs'}=\emptyset$
($\zw_{\zs}\cap\zw_{\zs'}\not=\emptyset$).
Note that if the scales of $\z8T$  differ by a factor of $2^{2j}$ [That is,
if $\abs{I_\zs}<\abs{I_{\zs'}}$, then $2^{2j}\abs{I_\zs}\le\abs{I_{\zs'}}$ for
all $\zs,\zs'$.]  then a tree $\z8T$ can be uniquely decomposed as a union of a
$1$--tree and a $2$--tree.   [We also remark that these definitions play a role
that is parallel to the notions of a tree  in \cite{ltcar}.]


\def\engs{ \text{\rm size}(\z5S)}    \def\eng#1]{\text{size}(#1)}


For $\z7S\subset \z6P$, define the {\em ``size of $\z7S$"} to be
\md0
\eng \z7S]:={}
\sup_{\z8T\subset \thin\z7S]}\Bigl[ \abs{I_{\z8T}}^{-1}\sum_{\zs\in\z8T}
\abs{\ip f,\zf_\zs.}^2\Bigr]^{1/2}
\emd
where the supremum is formed over all $1$--trees $\z8T\subset
\thin\z7S]:=\bigcup_{s\in\z5S}\thin \z5s] $.    The
central lemma concerning size is


\bthm l.eng   A finite collection $\z7S\subset\z6P$ is a union of
collections $\z7S(n)$, $n\in\ZZ$ for which
$\eng \z7S(n)]\le{}j2^n$ and
\md1  \Label e.eng
\sum_{ \z5s\in\z7S(n)^*}\abs{I_{\z5s}}\seq{}2^{-2n},
\emd
where $\z7S(n)^*$ consists of the maximal $\z5s\in \z7S(n)$.
\ethm

\noindent Observe that $j$ (that is a measure of how fat the tiles are) enters into
this lemma, albeit in a weak fashion.

\medskip

For $\z8S\subset \z9P$ set
\md0
M_j^{\z8S}=\sum_{\zs\in \z8S}\ip f,\zf_\zs.\zvf_\zs.
\emd
If $\z7S\subset \z6P$, define $M_j^{\z7S}$ to be $M_j^{\thin \z7S]}$.
Concerning trees, our central lemma is


\bthm  l.tree    For all trees $\z7T$,
\md5 \Label e.tree1
\norm M_j^{\z7T}.p.\seq\eng \z7T]
j2^{j(3-5/p)}\abs{I_{\z7T}}^{1/p},\qquad 2<p<\zI.
\\ \Label e.tree2
\abs{ M_j^{\z7T}(x)}\seq 2^{2j}\eng \z7T]\zq_{I_{\z7T}}(x)^m,\qquad x\not\in
2I_{\z7T}, \ m\ge1,
\emd
where
$\zq_J(x)=(1+\text{dist}(x,J)\abs{J}^{-1})^{-1}$.
\ethm
\noindent   Notice that  the first estimate should  be compared to Stein's
estimate for $\z0D[M_j]$, and is only slightly worse than that estimate if $p=2$.
That the (large) factor of $2^{2j}$ enters into the second estimate is completely
harmless.

Set $\zve=(200)^{-1}$, $p=9/4$, $\zm=7/9$ and $\zg=8/9$.

For
$n>-\zg j$, we in essence relie upon the fact that
$M_j^{\z7S(n)}$ is supported on a set of small measure.  To make
this precise,
let $E_n=\bigcup_{s\in \z7S(n)^*}2^{\zve j}I_{\z5s}$.
And set $F_0=\bigcup_{n>-\zg j}E_n$.  This set has measure
\md0
\abs{F_0}\le{}\sum_{n>-\zm j}\abs{E_n}\seq{}\sum_{n>-\zm j}2^{\zve
j-2 n}\seq 2^{2\zg j}.
\emd
We do not need to estimate $M_j$ on this set.  Using \e.tree2/ we see
that
\md0
\norm M_j^{\z7S(n)}.L^1(F_0^c).\seq 2^{-100j}
 \sum_{s\in \z7S(n)^*}\abs{I_{\z5s}}\seq{}
2^{-100 j-2n}.
\emd
Bringing these estimates together, we see that  for the collection
$\z7{\tilde S}=\bigcup_{n>-\zm j}\z7S(n)$,
we have  $\abs{\{M_j^{\z7{\tilde S}}>1\}}\seq2^{2\zg j}$, as is required in
\e.d2o/.

\bigskip

For $n\le-\zm j$, we need a more involved argument.  We encode some of
the necessary combinatorics into this Lemma.


\bthm l.s
For $n\le-\zm j$ there is a set $E_n\subset \ZR$ with $\abs{E_n}\seq
2^n$ so that the collection $\z7{\tilde S}(n)=\{\z5s\in \z7S(n)\mid
I_{\z5s}\not\subset E_n\}$ is a union of collections $\z7U(n,k)$,  $1\le k\le -500
n$, which satisfies these properties.  For each $1\le k\le -500n$,
\begin{itemize}
\item[(i)]  $\z7U(n,k)$ is uniquely decomposable into maximal disjoint
trees $\z7T\in\z6T({n,k})$.
\item[(ii)]  $\norm \sum_{\z7T\in\z6T({n,k})} 1_{I_{\z7T}}.\zI.\seq2^{-10n}$.
\item[(iii)]  $\{2^{-\zve n}I_{\z7T}\times \zw_{\z7T}\mid
\z7T\in\z6T({n,k})\}$ are pairwise disjoint rectangles.
\item[(iv)]  For all $\z5s\in\z7U(n,k)$
\md0
I_{\z5s}\not\subset\bigcup_{\z7T\in\z6T({n,k})}\{x\mid
\text{dist}(x,I_{\z7T})<2^{10n}\abs{ I_{\z7T}}\}.
\emd
\item[(v)]  Either $\z7T=\{I_{\z7T}\times \zw_{\z7T}\} $ for all $
\z7T\in\z6T({n,k}) $ or
$I_{\z7T}\times\zw_{\z7T}\not\in\z7T$ for all $\z7T\in\z6T({n,k})$.
\item[(vi)]
If $s,s'\in \z7U(n,k)\cup\bigcup\{I_{\z7T}\times\zw_{\z7T}\mid \z7T\in\z6T\}$  and
$\abs{I_{\z5s}}<\abs{I_{\z5s'}}$ then
$2^{-200n}\abs{I_{_{\z5s}}}\le\abs{I_{\z5s'}}$.
\end{itemize}
\ethm

We do not estimate $M_j^{\z5S}$ on the set $F_1=\bigcup_{n\le{}-\zm
j}E_n$. As this set has measure $\abs{F_1}\seq{}2^{-\zm
j}\seq{}2^{-2\zg j}$, there is no harm in doing this.

 Off of this set, our lemma permits the following construction.
For all $n\le-\zm j$, $0\le k\le-500n$ and $\z7T\in\z6T({n,k})$ there is
a functions $N^{\z7T}$ for which
\md5 \Label e.n1
\abs{N^{\z7T}(x)-M_j^{\z7T}(x)}\seq{}2^{10n}\zq_{I_{\z7T}}(x)^{100},\qquad
x\in\ZR,
\\  \Label e.n2
\text{ The functions $N^{\z7T}$ are disjointly supported in $\z7T\in\z6T({n,k})$.}
\emd
But then we can estimate by \e.tree1/
\md2
\Bigl\lVert \sum_{\z7T\in\z6T({n,k}) } N^{\z7T}\Bigr\rVert_p^p={}&
        \sum_{\z7T\in\z6T({n,k}) } \lVert N^{\z7T}\rVert_p^p
\\
{}\seq{}&  j^p2^{np+(3-5/p)pj} \sum_{\z7T\in\z6T({n,k})
}\abs{I_{\z7T}} \\  {}\seq{}&j^5 2^{7j/4+n/4}.
\emd
Thus certainly
\md0
\abs{\{M^{\z7U(n,k)}\ge2^{n/18}\}}\seq{}2^{7j/4+n/8}.
\emd
This is summable over  $n\le-\zm j$ and $0\le{}k\le{}-500n$ and so
completes our proof of \e.d2o/.

\smallskip

[This   interplay between $L^2$ and $L^p$
estimates is due to C.~Thiele \cite{t} and   contrasts with  the
argument of Lacey and Thiele \cite{ltcar}.
 The latter paper uses two notions  ``energy" (the current
``size") and ``mass", which are in some sense dual to one
another. The notion of ``mass"   seems to have little utility in
this paper: ``Mass" can be exploited through devices linked to
the Hardy--Littlewood maximal function, but our kernels bear no
close connection to that maximal function.]

\medskip

The construction relies on an argument from \cite{lmax}.  Fix $n,k$, set
$\z7U:=\z7U(n,k)$
and $\z6T:=\z6T({n,k})$.   To each $\z5s\in\z7U$ we construct a set
$G_{\z5s}$ as follows.  Recall
$(v)$ from the Lemma.  If each $\z7T\in\z6T$ consists only of a top we
set $G_{\z5s}={}2^{-\zve n}I_{\z5s}$ where $s$ is the top of the tree and
$N^{\z7T}= 1_{G_{\z5s}}\sum_{\zs\in\thin \z5s]}\ip f,\zf_\zs. \zvf_\zs$.  Then \e.n2/ follows from $(iii)$ and \e.n1/
follows from \e.tree2/.

We thus assume that no tree $\z7T\in\z6T$ contains its top. We then
make the following definitions for
$\z5s\in\z7T$.
\md4
G_{\z5s}=I_{\z7T}-\bigcup_{\z7T'\in\z6T(\z5s)}I_{\z7T'},
\\
\z6T(\z5s):=\{\z7T'\in\z6T-\{\z7T\}\mid  \zw_{\z5s1}\supset\zw_{\z7T'},\
I_{\z7T'}\subset I_{\z7T}\},
\\
N^{\z7T}=\sum_{\z5s\in\z7T} 1_{G_{\z5s}}\sum_{\zs\in\thin \z5s]}\ip f,\zf_\zs.
\zvf_\zs
\emd


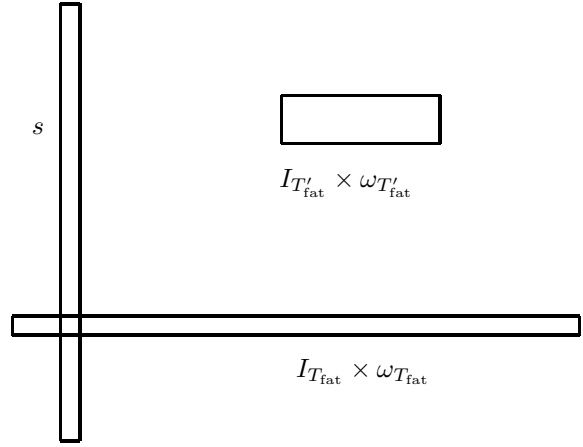
\begin{figure}[htbp]
\begin{center}
 \setlength{\unitlength}{0.00033333in}
\begingroup\makeatletter\ifx\SetFigFont\undefined%
\gdef\SetFigFont#1#2#3#4#5{%
  \reset@font\fontsize{#1}{#2pt}%
  \fontfamily{#3}\fontseries{#4}\fontshape{#5}%
  \selectfont}%
\fi\endgroup%
{\renewcommand{\dashlinestretch}{30}
\begin{picture}(8916,6906)(0,-10)
\thicklines
\drawline(33,1983)(8883,1983)(8883,1683)
	(33,1683)(33,1983)
\drawline(783,6858)(1083,6858)(1083,33)
	(783,33)(783,6858)
\drawline(4233,5433)(6708,5433)(6708,4683)
	(4233,4683)(4233,5433)
\put(333,4833){\makebox(0,0)[lb]{\smash{{{\SetFigFont{10}{12}{\rmdefault}{\mddefault}{\updefault}$s$}}}}}
\put(4458,1033){\makebox(0,0)[lb]{\smash{{{\SetFigFont{10}{12}{\rmdefault}{\mddefault}{\updefault}$I_{\z7T}\times
\zw_{\z7T}$}}}}}
\put(4200,4008){\makebox(0,0)[lb]{\smash{{{\SetFigFont{10}{12}{\rmdefault}{\mddefault}{\updefault}$I_{\z7T'}\times
\zw_{\z7T'}$}}}}}
\end{picture}
}\end{center}
\Label f.truncate
\caption{\small The top of the tree $\z7T'$  does not intersect the set $G_s$
since the intervals $I_{\z7T}$ and $I_{\z7T'}$ intersect.}
\end{figure}

We verify \e.n2/.
Since the support of $\zvf_\zs$ is in $\{x\mid N(x)\in\zw_{s1}\}$, \e.n2/
is a consequence of the observation that if $G_{\z5s}\times
\zw_{\z5s1}\cap G_{\z5s'}\times \zw_{\z5s'1}\not=\emptyset$ then $\z5s$ and
$\z5s'$ are in the same tree.  Indeed write $\z5s\in\z7T$ and
$\z5s'\in\z7T'$ and assume say $\zw_{\z5s'1}\subset\zw_{\z5s1}$.
If $\zw_{\z5s'1}=\zw_{\z5s1}$ and the two trees are distinct then $I_{\z7T}$ and
$I_{\z7T'}$ are disjoint by $(i)$.  Assume $\zw_{\z5s'1}\subset_{\not=}\zw_{\z5s1}$ and
$G_{\z5s}\cap G_{\z5s'}\not=\emptyset$.
Then
$\zw_{\z7T},\zw_{\z7T'}\subset \zw_{\z5s1}$. See figure 2. The tops $I_{\z7T}$ and $I_{\z7T'}$ must
intersect.  Assuming $I_{\z7T}\subset I_{\z7T'}$ then
$\z5s<I_{\z7T'}\times \zw_{\z7T'}$.  But then condition $(i)$ forces $\z7T=\z7T'$.  Thus we
must have $I_{\z7T'}\subset I_{\z7T}$, which by definition means that
$I_{\z7T'}\cap G_{\z5s}=\emptyset$, so that $G_{\z5s}\cap
G_{\z5s'}=\emptyset$.  This is a contradiction and so proves \e.n2/.

We verify \e.n1/.  In the case of $x\not\in I_{\z7T}$ this follows from
\e.tree2/ and conditions $(iv)$ and $(vi)$ of \l.s/.  We do not comment
further.  For $x\in I_{\z7T}$ we in fact have
$N^{\z7T}(x)=M_j^{\z7T}(x)$ unless $x\in I_{\z7T'}$
and $\z7T'\in\z8T(s)$ for some $s\in\z7T$.
Indeed, with $\z7T$ fixed we can assume that $I_{\z7T'}\subset I_{\z7T}$ for
all $\z7T'\in\z6T$.  Then
we shall  just reverse the order of summation below.
\md0
\abs{N^{\z7T}(x)-M_j^{\z7T}(x)}\le{}\sum_{\z7T'\in\z6T-\{\z7T\}}
\sum_{\zs\in\thin \z6T(\z7T')]}
  \abs{\ip f,\zf_\zs.\zvf_\zs(x)},
\emd
where $\z6T(\z7T'):=\{\z5s\in\z7T\mid \z7T'\in\z6T(\z5s)\}$ and
$\z6T(\z5s)$ was used to define
$G_{\z5s}$.
But again condition $(iv)$ and $(vi)$ imply that
\md0
\sum_{\zs\in\thin \z6T(\z7T')]} \abs{ \ip
f,\zf_\zs.\zvf_\zs(x)}\seq2^{200n}1_{\z7T'}(x)
\emd
and \e.n1/ follows from condition $(ii)$.

\medskip
Our proof of \e.d2o/ is complete modulo the proofs of the lemmas,
which are taken up in the next section.


\section{Proofs of the Lemmas}

\begin{proof}[Proof of \l.eng/]
The argument is a variant of one in \cite{icm} and has been used several
times since.  We give the details, although only small changes are needed to account
for the disparity between fat and thin tiles.   The most expedient treatment requires   a new
definition of a tree.

Fix a choice of integer $0\le{}k<200j$.
For a $1$--tree $\z8T$ call a subset $\z8T^{\ell}\subset \z8T$ a
{\em
left--tree (right--tree) }if there is a
$\zx_{\z8T}\in\zw_{\z8T}$ with $\zx_{\z7T}$ to the left (right)
of every $\zw_s$, $s\in\z8T^\ell$. In addition require that for
all $s\in\z8T^\ell\cup\{I_{\z8T^\ell}\times \zw_{\z8T^\ell}\}
$, $\log_2\abs{I_s}\in k+200j\ZZ$. Define
``the left size of $\z7S$," or
$\text{$\ell$--size}(\z7S)$ as
\md0
\text{$\ell$--size}
(\z7S)=\sup\Bigl\{\Bigl[\abs{I_{\z7T}}^{-1}\sum_{s\in\z7T^\ell}
\abs{\ip f,\zf_\zs.}^2\Bigr]^{1/2}\Bigr\}
\emd
where the supremum is over all left--trees $\z7T^\ell$ with
$\z7T\subset \z7S$.

We prove this statement.  For any finite $\z7S\subset \z6P$ set
$\zve=\text{$\ell$--size} (\z5S)$.  Then $\z7S=\z5S_{\text{lo}}\cup
 \z5S_{\text{hi}}$ with $\text{$\ell$--size}
(\z5S_{\text{lo}})\le\zve/4$ and $\z5S_{\text{hi}}$ is a union of trees $\z7T\subset \z6T$ with
\md1 \Label e.e1
\sum_{\z7T\in\z6T}\abs{I_{\z7T}}\seq\zve^{-2}.
\emd

An inductive application of this statement proves \l.eng/ with $\eng \z7S] $
replaced by $\text{$\ell$--size}(\z7S)$.
The factor $j$ does not enter into this statement of the lemma.
 The same statement is
true for right--size.  Letting $k$ vary from $0$ to $200j$ proves
the Lemma as stated.

\medskip

The construction of $\z5S_{\text{hi}}$ and $\z6T$ is inductive.  The construction
also associates to each $\z7T\in\z6T$ a particular left--tree
$\z8T^\ell$ which are used to prove \e.e1/.  Initially set
$\z5S^{\text{stock}}:=\z7S$.
Select a tree $\z7T\subset \z5S^{\text{stock}}$ so that
\begin{itemize}
\item[$(a)$]  $\z7T$ contains a left--tree $\z8T^\ell$ with
\md0
\sum_{\zs\in\z8T^\ell}\abs{\ip f,\zf_\zs.}^2\ge
(4\zve)^{-2} \abs{I_{\z7T}}.
\emd

\item[$(b)$]  $I_{\z7T}$ is maximal amoung trees satisfying condition  $(a)$
and $\z7T$ is the maximal tree in $\z5S^{\text{stock}}$ with that top.

\item[$(c)$] $\zx_{\z8T^\ell}$ is right--most amoung trees satisfying $(a)$ and $(b)$.
\end{itemize}
Then add $\z7T$ to $\z6T$, set
$\z5S^{\text{stock}}:=\z5S^{\text{stock}}-\z7T$.  Repeat this procedure until there is no tree satisfying
$(a)$.  Then set $\z5S_{\text{lo}}:=\z5S^{\text{stock}}$.  By definition, $\ell-\eng
\z5S_{\text{lo}}] \le\zve/4$.

The left--trees we have constructed satisfy this disjointness property.
For $\z7T\not=  \z7T'\in\z6T$ and $\zs\in\z8T^\ell$ and $\zs'\in\z8T^{\ell'}$
\md1\Label e.e2
\text{ if $\zw_\zs\subset_{\not=}\zw_{\zs'}$ then $I_{\z7T}\cap
I_{\zs'}=\emptyset$.}
\emd
Indeed, $\zw_{\z7T}\subset\zw_{\zs'}$ so that $\zx_{\z7T'}<\zx_{\z7T}$.
Thus the tree $\z7T$ was constructed before $\z7T'$.  But
if $I_{\z7T}\cap I_{s'}\not=\emptyset$ we see that
${s'}<I_{\z7T}\times \zw_{\z7T}$ where $\zs\in\thin s']$.  Hence $s'\in\z7T$ which is a
contradiction.  See figure 3.


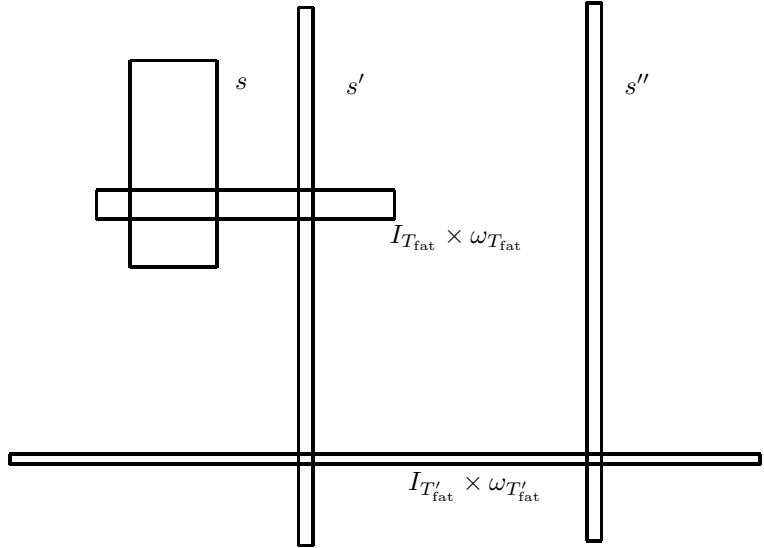
\begin{figure}[htbp]
\begin{center}
 \setlength{\unitlength}{0.00033333in}
\begingroup\makeatletter\ifx\SetFigFont\undefined%
\gdef\SetFigFont#1#2#3#4#5{%
  \reset@font\fontsize{#1}{#2pt}%
  \fontfamily{#3}\fontseries{#4}\fontshape{#5}%
  \selectfont}%
\fi\endgroup%
{\renewcommand{\dashlinestretch}{30}
\begin{picture}(11766,8556)(0,-10)
\thicklines
\drawline(1383,5583)(6033,5583)(6033,5133)
	(1383,5133)(1383,5583)
\drawline(33,1458)(11733,1458)(11733,1308)
	(33,1308)(33,1458)
\drawline(1908,7608)(3258,7608)(3258,4383)
	(1908,4383)(1908,7608)
\drawline(4758,8433)(4533,8433)(4533,33)
	(4758,33)(4758,8433)
\drawline(9258,8508)(9033,8508)(9033,108)
	(9258,108)(9258,8508)
\put(5658,4758){\makebox(0,0)[lb]{\smash{{{
\SetFigFont{10}{10}{\rmdefault}{\mddefault}{\updefault}
$I_{\z7T}\times \zw_{\z7T}$}}}}}
\put(6108,900)
{\makebox(0,0)[lb]{\smash{{{\SetFigFont{10}{12}{\rmdefault}{\mddefault}{\updefault}
$I_{\z7T'}\times \zw_{\z7T'}$ }}}}}
\put(3408,7158){\makebox(0,0)[lb]{\smash{{{\SetFigFont{10}{12}{\rmdefault}{\mddefault}{\updefault}
$s$}}}}}
\put(5133,7083){\makebox(0,0)[lb]{\smash{{{\SetFigFont{10}{12}{\rmdefault}{\mddefault}{\updefault}
$s'$}}}}}
\put(9483,7083){\makebox(0,0)[lb]{\smash{{{\SetFigFont{10}{12}{\rmdefault}{\mddefault}{\updefault}
$s''$}}}}}
\end{picture}
}
\end{center}
\Label  f.disjointtrees
\caption{ By the manner in which the trees are constructed, the tree $\z7T$
was constructed before $\z7T'$.  Hence the tile $s'$ must be in $\z7T$.  But the
tile $s''$ is a member $\z7T'$.}
\end{figure}

\def\redt{T^{\text{red}}}

Let $\redt$ be those $\zs\in\z8T^\ell$ for which if
$\abs{I_\zs}<\abs{I_{\redt}}$ then $\text{dist}(I_\zs,\partial
I_{\redt})\ge\frac1{32}\abs{ I_{\redt}}$. [``red" is for ``reduced." 
Note that the top is
permitted to be in $\redt$.  And that if
$\abs{I_\zs}<\abs{I_{\redt}}$ then $\abs{I_\zs}$ is in fact much
smaller than $\abs{I_{\redt}}$.]  As $\ell$--$\eng
\z7S]\le{}\zve$, it follows that \md0 \sum_{\zs\in\redt}\abs{\ip
f,\zf_\zs.}^2\ge{}\frac{\zve^2}{32}\abs{I_\redt}. \emd

Set $\z8S=\bigcup_{\z7T\in\z6T}\redt $. And
\md0
B:=\NOrm \sum_{\zs\in\z8S}{\ip f, \zf_\zs.} \zf_\zs.2.
\emd
Observe that by Cauchy--Schwartz and $\norm f.2.=1$,
\md3  \nonumber
\frac {\zve^2}{32}\sum_{\z7T\in\z6T}\abs{I_{\z7T}}\le{}&  \sum_{\zs\in\z8S}
\abs{\ip f,\zf_\zs.}^2
 \\ \nonumber
{}={}& \IP f,{ \sum_{\zs\in\z8S}\ip f,\zf_\zs. \zf_\zs } .
\\
\Label e.e3   {}\le{}&B.
\emd
To conclude \e.e1/ we show that
\md1 \Label e.e4
B^2\seq{} \zve^2\sum_{\z7T\in\z6T}\abs{I_{\z7T}}.
\emd

By expanding the $L^2$ norm $B^2\le{}2(B_1^2+B_2^2)$ where we define
\md4
B_j^2:=\sum_{\zs\in\z8S}\ip f,\zf_\zs. \sum_{\zs'\in\z8S^j(\zs)}\ip
\zf_\zs,\zf_{\zs'}.\ip \zf_{\zs'},f.
\\
\z8S^1(\zs):=\{\zs'\in\z8S\mid \zw_\zs=\zw_{\zs'}\},\qquad
    \z8S^2(\zs):=\{\zs'\in\z8S\mid \zw_\zs\subset_{\not=}\zw_{\zs'}\}.
\emd
Note that if $\zw_\zs\subset\zw_{\zs'}$ we have
\md1\Label e.e5
\abs{\ip  \zf_\zs,\zf_{\zs'}.}\seq \sqrt{\frac {\abs {I_{\zs'}}}{\abs{I_\zs}}  }
\zq_{I_\zs}(c(I_{\zs'}))^{100}.
\emd

\medskip

To bound $B_1$ fix a dyadic interval $\zw$. This last estimate
and Cauchy--Schwartz estimate shows that
\md0
\ABs{  \sum_{ \substack{ \zs\in\z8S \\ \zw_\zs=\zw}}
        \ip f,\zf_\zs. \sum_{\zs'\in\z8S^1(s)}
	\ip \zf_\zs,\zf_{\zs'}.\ip \zf_{\zs'},f.
}\seq{}
 \sum_{ \substack{ \zs\in\z8S \\ \zw_\zs=\zw}} \abs{ \ip f,\zf_\zs.}^2.
\emd
Hence by \e.e3/ and summing over $\zw$,
\md0
B_1^2\seq\sum_{\zs\in\z8S}\abs{ \ip f,\zf_\zs.}^2
{}\seq{}  \zve^2\sum_{\z7T\in\z6T}\abs{I_{\z7T}}.
\emd
as $\text{$\ell$--size}(\z8S)=\zve$.  This is the first step in establishing \e.e4/.

To control $B_2^2$ we must use the disjointness property \e.e2/.  Fix a
tree $\redt$ and consider $\zs\in\redt$.  Then the intervals
$\{I_{\zs'}\mid \zs'\in\z8S(s)\}$ are pairwise disjoint and contained in
$(I_{\z7T})^c$.  To see this note that for all $\zs',\zs''\in\z8S(s)$ we have
$\zw_\zs\subset \zw_{\zs'}\cap \zw_{\zs''}$.  So \e.e2/ implies $I_{\zs'}\cap
I_{\zs''}=\emptyset$.  Then we can estimate
\md2
\sum_{\zs'\in\z8S(\zs)}\abs{  \ip f,\zf_\zs.\ip \zf_\zs,\zf_{\zs'}.\ip
\zf_{\zs'},f.  }
{}\seq{}& \zve^2 \sum_{\zs'\in\z8S(s)}
\zq_{I_\zs}(c(I_{\zs'}))^{100}\abs{I_{\zs'}}
\\ {}\seq{}&
 \zve^2\int_{ (I_{\z7T})^c} \zq_{I_\zs}(x)^{90}\; dx
 \\{}\seq{}&
 \zve^2 \Bigl(\frac{\abs{I_\zs}}{\abs{I_{\redt}}}\Bigr)^{10}\abs{I_\zs}.
\emd
Here, we have in addition relied upon the estimate $\abs{ \ip
f,\zf_\zs.}\le\zve\sqrt{\abs{I_\zs}}$.  Finally, the estimate
below follows as $I_\zs$ is both much smaller than $I_\redt$ and
not close to the boundary of $I_\redt$.  This  completes the
proof of \e.e4/. \md0
\sum_{\zs\in\z8T^\ell} \int_{(I_{\z7T})^c} \zq_{I_\zs}(x)^{90}\; dx \seq{}
\abs{I_{\z7T}}.
\emd

\end{proof}


\begin{proof}[Proof of \l.tree/.]
We begin by verifying \e.tree2/.  For any $\zs\in\z9P$ and $m\ge0$
observe that
\md1\Label e.1+
\abs{\zvf_\zs(x)}\seq{}  \abs{I_\zs}^{-1/2}\zq_{I_\zs}(x)^m,\qquad x\in\ZR.
\emd
Indeed, after taking dilation and translation into account this
estimate reduces to
\md0
\ABs{  2^{j/2}\int \tilde\zc (2^j\zx+\zx_0)\ex \zx^2+\zx y]\; d\zx
}\seq{}
2^{-j/2}(1+2^j\abs y)^{-m},\qquad y\in\ZR.
\emd
Here, $\tilde\zc$ is a Schwarz function supported in
$\frac12\le\abs\zx\le2$ and $2^{j-1}\le
\zx_0\le2^{j+1}$.   But then at most ${}\seq1$ oscillations  of $\ex
\zx^2]$ are relevant to the integral, so the estimate follows by a repeated
integration by parts.  Then   \e.1+/ plus a routine argument proves \e.tree2/.

\medskip

Turning to the estimate \e.tree1/,
note that any tree $\z8T\subset\z9P$ is a union of ${}\seq{}j$ $1$ and $2$---trees.
It suffices to prove \e.tree1/ without the leading factor of $j$ on the right for
$1$ and $2$---trees.

We consider first the case of a
$2$--tree $\z8T$.  In this
case, the sets $\zw_{\z5s1}$ for $\z5s\in\z8T$ are disjoint and for
$\zs\not=\zs'\in\z8T $, we have either $\abs{I_\zs}\not=\abs{I_{\zs'}}$, in which
case $\zvf_\zs$ and $\zvf_{\zs'}$ are disjointly supported,  or $I_\zs\cap
I_{\zs'}=\emptyset$, in which case we relie upon the decay \e.1+/.   Thus,
\md0
\abs{ M_j^{\z8T }(x)}\seq{} \eng \z8T]\zq_{I_{\z8T}}(x)^8.
\emd
That is, \e.tree1/ is trivially satisfied in this case. [This argument
is the key motivation for the definitions of $\z6P$ and $\z9P$.]

\smallskip

We now turn to the case of a $1$---tree $\z8T$.  A specific case unlocks the general case.
Suppose that $\z8T$ is a tree with  $\abs{I_\zs}=\abs{I_{\zs'}}$ for all
$\zs,\zs' \in \z8T$ and $0=c(\zw_{\zs,1})$.
  Then from \e.dj/ we have for all
$2<p<\zI$, \md3   \nonumber
\norm M^{\z8T}.p.\le{}&  \Norm D_j \Bigl(\sum_{s\in\z8T} {\ip f,\zf_\zs.} \zf_\zs
\Bigr).p.
\\  \nonumber {}\seq{}& 2^{j(1-1/p)}\Norm \sum_{\zs\in\z8T} {\ip
f,\zf_\zs.} \zf_\zs .p. \\ \Label e.pt2
{}\seq{}& 2^{j(1-1/p)+2j(1-2/p)}\eng \z8T]  \abs{I_{\z8T}}^{1/p}.
\emd
The last line follows as there are $2^{2j}$ tiles $\zs$ in any $\thin s]$ for
$\z5s\in\z6P$.  Combine this with a trivial interpolation
argument to conclude this case.

\medskip

More generally, for any tree $\z8T$, observe that there is a connection
to the space of functions of bounded mean oscillation.  The distinction 
between fat and thin tiles must enter into this relationship however. 
And in particular it is 
\md1
\Label e.pt1
\Norm \sum_{\zs\in\z8T} {\ip f,\zf_\zs.} \zf_\zs
.BMO.\seq2^{2j}\eng \z8T].
\emd
This follows from the definition of size.

Let $\zw_{\z8T}=\zw^1\subset_{\not=}\zw^2\subset_{\not=}\cdots$ be the 
maximal sequence of dyadic intervals containing $\zw_{\z8T}$.   Let 
$\z8T^l=\{  
\z5s\in\z8T\mid \zw_\z5s=\zw^l\}$.  Then the functions $M^{\z8T^l}$ are 
disjointly supported in $l$.  
 Hence, 
\md2
\lVert M^{\z8T}\rVert_p^p={}&  \sum_l\lVert M^{\z8T^l}
\rVert_p^p
\\ {}\seq{}   &
2^{jp(1-1/p)} \sum_l  \Bigl\lVert \sum_{\zs\in\z8T^l} \ip
f,\zf_\zs. \zf_\zs \Bigr\rVert_p^p
\\  {}\seq{}   &
2^{jp(1-1/p)}  \Bigl\lVert \Bigl[ \sum_l  \ABs{
\sum_{\zs\in\z8T^{l}} \ip f,\zf_\zs. \zf_\zs }^2 \Bigr]^{1/2}
\Bigr\rVert_p^p \\   {}\seq{}   &2^{jp(3-5/p)}\eng
\z8T]\abs{I_{\z8T}}^{1/p}. \emd
Here we relie on $p>2$, \e.pt1/ and\e.pt2/.  
\end{proof}

\begin{proof}[Proof of \l.s/.]
 We shall show that
there is a set $E_n\subset \ZR$ so that
$\abs{E_n}\seq{}2^n$ and the collection 
 $\z7{\widetilde S}(n):=\{s\in\z7S(n)\mid I_s\not\subset E_n\}$
  is a union of collections
$\widetilde{\z7U}(n,k)$, $0\le k\le{}-50n$ satisfying  $(i)$---$(iii)$.

The last three conditions of the Lemma are trivially satisfied by making further
subdivisions of the   subcollections $\z7U(n,k)$, and making a small
further contribution to the exceptional set $E_n$.  Thus, the
lemma will follow in complete generality.

\medskip

Fix $n$ and set $\z7S=\z7S(n)$.  Condition $(ii)$ is also easy to satisfy.  For
the first contribution to our exceptional set, define 
\md0
E^1:=\Bigl\{x\mid \sum_{x\in\z7S^*}1_{I_s}(x)>2^{-10n}\Bigr\}
\emd
where $\z7S^*$ consists of the maximal elements of $\z7S$.  By \e.eng/ $\abs {E^1}
\seq2^n$.   
We can assume that for all $s\in\z7S$, $I_s\not\subset E^1$.  Then certainly
$(ii)$ is true.

We now show that $\z7S$ is decomposable into subcollections $\z7U(k)$,
$1\le{}k\le{}-10n$ which are uniquely decomposable into maximal trees.  This last
condition is true iff to each $s\in\z7U(k)$ there is a unique
maximal $s^*\in\z7U^{k}$ with $s<s^*$.  And this is so iff the
collections $\z7U(k)$ does not admit a vee in the partial order
on tiles.  A {\em vee} is three tiles $s,s',s''$ with $s<s',s''$
but $s'$ and $s''$ are not comparable with respect to the partial
order on tiles. 

To acheive this, we employ a method of Fefferman \cite{f}.
Define a counting function \md0
C(s):=\sharp \{s^*\in\z7S^*\mid s<s^*\}.
\emd
Then $C(s)\le{}2^{-10n}$ for all $s\in\z7S$ as $(ii)$ is true.  Take the sets $\z7U(k)$ to be
$\{s\mid 2^{k-1}\le{}C(s)<2^k\}$.

That these sets do not contain vees follows immediately from the observation 
 that $C(s)$ is superadditive in
this sense.  If $s,s',s''\in\z7S$ is a vee, then $C(s)\ge{}C(s')+C(s'')$.   
[Then if $s',s''\in\z7U(k)$ we see that
$C(s)\ge2^{k-1}+2^{k-1}=2^k$, so it can not be in $\z7U(k)$.]
Indeed, there there can be no  maximal tile $s'''$ larger than both $s'$ and $s''$, for this
would force $s'$ and $s''$ to be comparable in the partial order, as one checks
immediately.  Hence the maximal tiles greater than $s'$ are disjoint from those
greater than $s''$, which proves the superadditivity property. 

The last condition to verify is $(iii)$, which requires another class of
contributions to the exceptional set.  Fix a choice of $1\le{}k\le{}-10n$. 
 Consider the maximal tiles $\z7U^*(k)$.  We want to separate
these tiles after expanding the coordinates $I_s$ by a factor of $2^{-\zve n}$.  
This can be done, up to an exceptional set and a further division of $\z7U^*(k)$,
 by applying Lemma~$4.4$ to $S=\z7U^*$,
with $A=2^{\zve n}$.   The details are omitted. 

\end{proof}

   \endgroup  

\bigskip 
 {\parindent=0pt\baselineskip=12pt\obeylines 
Michael T. Lacey 
School of Mathematics 
Georgia Institute of Technology 
Atlanta GA 30332
\smallskip
\tt lacey@math.gatech.edu
\tt http://www.math.gatech.edu/\~{}lacey }

\end{document}